\documentstyle[amssymb,12pt]{article}

\input{tcilatex}

\def\text{\mbox}

\begin{document}

\begin{center}
\mathstrut

{\Large {\bf Notes on Matrix Valued Paraproducts}}

\medskip

Tao MEI \footnote{%
Partially supported by the National Science Foundation 0200690.}
\end{center}

\begin{quotation}
{\bf Abstract } Denote by $M_n$ the algebra of $n\times n$ matrices. We
consider the dyadic paraproducts $\pi _b$ associated with $M_n$ valued
functions $b$, and show that the $L^\infty (M_n)$ norm of $b$ does not
dominate $||\pi _b||_{L^2(\ell _n^2)\rightarrow L^2(\ell _n^2)}$ uniformly
over $n$. We also consider paraproducts associated with noncommutative
martingales and prove that their boundedness on bounded noncommutative $L^p-$%
martingale spaces implies their boundedness on bounded noncommutative $L^q-$%
martingale spaces for all $1<p<q<\infty $.
\end{quotation}

\section{Introduction}

\setcounter{theorem}{0}\setcounter{equation}{0}
Denote by $M_n$ the
algebra of $n\times n$ matrices. Let $({\Bbb {T},{\cal F}}_k,dt)$ be the
unit circle with Haar measure and the usual dyadic filtration. Let $b$
be an $M_n$ valued function on ${\Bbb {T}}$. The matrix valued dyadic
paraproduct associated with $b,$ denoted by $\pi _b,$ is the operator
defined as
\begin{equation}
\pi _b(f)=\sum_k(d_kb)(E_{k-1}f),\quad \forall f\in L^2(\ell _n^2).  \label{def}
\end{equation}
Here $E_kf$ is the conditional expectation of $f$ with respect to
${\cal F}_k,$ i.e. the unique ${\cal F}_k$-measurable function such
that
  $$
  \int_FE_kfdt=\int_Ffdt, \ \ \ \forall F\in {\cal F}_k.
  $$
And $d_kb$ is defined to be $E_kb-E_{k-1}b.$

In the classical case (when $b$ is a scalar valued function), paraproducts are usually
considered as dyadic singular integrals and play
important roles in the proof of the classical T(1) theorem. It is well
known that
\[
\left\| \pi _b\right\| _{L^2\rightarrow L^2}\backsimeq \left\|
b\right\| _{BMO_d},
\]
where BMO$_d$ denotes the dyadic BMO norm defined as
\[
\left\| b\right\| _{BMO_d}=\sup_m\|E_m\sum_{k=m}^\infty |d_kb|^2|\|
_{L^\infty}^{\frac12}.
\]
And by the Calder\'on-Zygmund decomposition and the Marcinkiewicz
interpolation theorem, we have $||\pi_b||_{L^p\rightarrow
L^p}\backsimeq ||\pi_b||_{L^p\rightarrow L^p}\backsimeq||b||_{BMO_d}
$ for all $1<p<\infty$.

When $b$ is $M_n$ valued, it was proved by Katz (\cite{K}) and
independently by Nazarov, Treil and Volberg (\cite{NTV}, see
\cite{[P]} for another proof by Pisier) that
\begin{equation}
\left\| \pi _b\right\| _{L^2(\ell _n^2)\rightarrow L^2(\ell _n^2)}\leq c\log
(n+1)\left\| b\right\| _{BMO_c}.  \label{katz}
\end{equation}
Here $\Vert \cdot \Vert _{{\rm BMO}_c}$ is the column BMO norm defined by
\[
\left\| b\right\| _{BMO_c}=\sup_m\left\| E_m\sum_{k=m}^\infty
(d_kb)^*(d_kb)\right\| _{L^\infty (M_n)}^{\frac 12},
\]
where $(d_kb)^{*}$ is the adjoint of $d_kb.$ Nazarov, Pisier, Treil and
Volberg (\cite{NPTV}) proved later that the constant $c\log (n+1)$ in (\ref
{katz}) is optimal. Thus the BMO$_c$ norm does not dominate $\left\| \pi
_b\right\| _{L^2(\ell _n^2)\rightarrow L^2(\ell _n^2)}$ uniformly over $n$.

Can we expect something weaker? In particular, does there exist a constant $c
$ independent of $n$ such that, for every $n\in {\Bbb N},$%
\begin{equation}
\left\| \pi _b\right\| _{L^2(\ell _n^2)\rightarrow L^2(\ell _n^2)}\leq
c\left\| b\right\| _{L^\infty (M_n)}?  \label{con}
\end{equation}
Some known facts made (\ref{con}) look hopeful. For example, the Hankel
operator associated with the $M_n$ valued function $b$ has a norm equivalent to $%
||b||_{(H^1(S^1))^{*}}$. Here $||\cdot ||_{(H^1(S^1))^{*}}$ denotes the dual
norm on the trace class valued Hardy space $H^1(S^1).$ And S. Petermichl
proved a close relation between $\pi _b$ and the Hankel operators associated
with $b$ $($see \cite{Pe}).

In this paper, we prove the following theorem, which shows there does not
exist any constant $c$ independent of $n$ such that (\ref{con}) holds.

\begin{theorem}
For every $n\in {\Bbb N},$ there exists an $M_n$ valued function $b$ with $%
\left\| b\right\| _{L^\infty (M_n)}\leq 1$ but such that
\[
\left\| \pi _b\right\| _{L^2(\ell _n^2)\rightarrow L^2(\ell _n^2)}\geq c\log
(n+1),
\]
where $c>0$ is independent of $n.$
\end{theorem}

This also gives a new proof that the constant $c\log (n+1)$ in (\ref{katz})
is optimal.

Denote by $S^p$ the Schatten $p$ class on $\ell ^2$.
For $f\in L^p(S^p),$ we define $%
\pi _b(f)$ as in (\ref{def}) also.
As pointed out in \cite{[P]}, it is easy to check that $\left\| \pi
_b\right\| _{L^2(S^2)\rightarrow L^2(S^2)}$ $=\left\| \pi _b\right\|
_{L^2(\ell ^2)\rightarrow L^2(\ell ^2)}$. For scalar valued $b$, as
we mentioned previously, we have $||\pi_b||_{L^p\rightarrow
L^p}\backsimeq||\pi_b||_{L^q\rightarrow L^q}.$ We wonder if this is
still true for matrix valued $b$, i.e. if $\pi_b$'s boundedness on
$L^p(S^p)$ implies their boundedness on $L^q(S^q)$ for all
$1<p,q<\infty$.

More generally, we can consider paraproducts associated with
noncommutative martingales. Let ${\cal M}$ be a finite von Neumann
algebra with a normalized faithful trace $\tau .$ For $1\leq p <
\infty$, we denote by $L^p({\cal M})$ the noncommutative $L^p$ space
associated with $({\cal M}, \tau)$. Recall the norm in $L^p({\cal
M})$ is defined as
\[
\left\| f\right\| _p=(\tau |x|^p)^{\frac 1p},\quad \forall f\in L^p({\cal M}%
),
\]
where $|f|=(f^{*}f)^{\frac 12}.$\ For convenience, we usually set $L^\infty (%
{\cal M})={\cal M\,}$equipped with the operator norm $\left\| \cdot \right\|
_{{\cal M}}.$ Let ${\cal M}_k$ be an increasing filtration of von Neumann
subalgebras of ${\cal M}$ such that $\cup _{k\geq 0}{\cal M}_k$ generates $%
{\cal M}$ in the w$^*-$ topology. Denote by $E_k$ the conditional expectation
of ${\cal M}$ with respect to ${\cal M}_k.$ $E_k$ is a norm 1 projection of $%
L^p({\cal M})$ onto $L^p({\cal M}_k)$.
For $1\leq p\leq \infty ,$ a sequence
$f=(f_k)_{k\geq 0}$ with $f_k\in L^p({\cal M}_k)$ is called a bounded
noncommutative $L^p$-martingale, denoted by $(f_k)_{k\geq 0}\in L^p({\cal M}%
),$ if $E_kf_m=f_k,\forall k\leq m$ and
\[
||(f_k)_{k\geq 0}||_{L^p({\cal M})}=\sup_k||f_k||_{L^p({\cal M})}<\infty .
\]
Because of the uniform convexity of the space $L^p({\cal M})$, for $%
1<p<\infty$, we can and will identify the space of all bounded $L^p({\cal M})
$-martingales with $L^p({\cal M})$ itself. In particular, for any $f\in L^p(%
{\cal M})$, set $f_k=E_kf$, then $f=(f_k)_{k\geq 0}$ is a bounded ${L^p(%
{\cal M})}$-martingale and $||(f_k)_{k\geq 0}||_{L^p({\cal M})}=||f||_{L^p(%
{\cal M})}$. Denote by
$d_kf=E_kf-E_{k-1}f.$

We say an increasing filtration ${\cal M}_k$ is ``regular" if there exists a
constant $c>0$ such that, for any $m, a\in {\cal M}_m,
a\geq 0$,
$$||a||_\infty\leq c||E_{m-1}a||_\infty.$$
For ${\cal M}$ with a regular filtration ${\cal M}_k$, $b\in L^2(\cal{M})$,
we define paraproducts $\pi _b,\widetilde{\pi }_b$ as operators for bounded ${L^p({\cal M})}$
$(1<p<\infty )$-martingales $f=(f_k)_{k\geq 0}$ as
\[
\pi _b(f)=\sum_kd_kbf_{k-1},\mbox{ \qquad }\widetilde{\pi }
_b(f)=\sum_kf_{k-1}d_kb.
\]

We prove the following result for $\pi _b$ and $\widetilde{\pi }_b:$

\begin{theorem}
Let $1<p<q<\infty $, if $\widetilde{\pi }_b$ and $%
\pi _b$ are both bounded on $L^p({\cal M})$ then they are both bounded on $%
L^q({\cal M}).$
\end{theorem}

We still do not know what happens when $p>q$.

\section{Proof of Theorem 1.1 and Application to ``Sweep" functions.}

Denote by $tr$ the usual trace on $M_n$ and $S_n^p (1\leq p<\infty$) the
Schatten $p$ classes on $\ell^2_n$.

\noindent{\bf Proof of Theorem 1.1.} Let $c(n)$ be the best constant such
that
\[
\left\| \pi _b\right\| _{L^2(\ell _n^2)\rightarrow L^2(\ell _n^2)}\leq
c(n)\left\| b\right\| _{L^\infty (M_n)},\quad \forall b\in L^\infty (M_n).
\]
Denote by $T$ the triangle projection on $S_n^1,$ we are going to show
\[
\left\| T\right\| _{S_n^1\rightarrow S_n^1}\leq c(n).
\]
Once this is proved, we are done since $\left\| T\right\| _{S_n^1\rightarrow
S_n^1}\backsim \log (n+1)$ (see \cite{KP}). Note that every $A$ in the unit ball of $S_n^1$
can be written as
\[
A=\sum_m\lambda ^{(m)}\alpha ^{(m)}\otimes \beta ^{(m)}
\]
with $\sum_m\lambda ^{(m)}\leq 1$, $\sup_m\{||\alpha ^{(m)}||_{\ell
_n^2},||\beta ^{(m)}||_{\ell _n^2}\}\leq 1.$ Therefore, we only need to show
\begin{equation}
\left\| T(\alpha \otimes \beta )\right\| _{S_n^1}\leq c(n)\left\| \alpha
\right\| _{\ell _n^2}\left\| \beta \right\| _{\ell _n^2},\ \ \forall \alpha
=(\alpha _k)_k,\beta =(\beta _k)_k\in \ell _n^2.  \label{T}
\end{equation}
Let $D$ be the diagonal $M_n$ valued function defined as
\[
D=\sum_{i=1}^nr_ie_i\otimes e_i
\]
where $r_i$ is the $i$-th Rademacher function on ${\Bbb {T}}$ and $(e_i)_{i=1}^n$
is the canonical basis of $\ell _n^2$. Given $\alpha =(\alpha _k)_k,\beta
=(\beta _k)_k\in \ell _n^2$, let

\[
f=D\alpha ,g=D\beta .
\]
Then $f,g\in L^2(\ell _n^2),$ and
\begin{equation}
\left\| f\right\| _{L^2(\ell _n^2)}=\left\| \alpha \right\| _{\ell
_n^2},\left\| g\right\| _{L^2(\ell _n^2)}=\left\| \beta \right\| _{\ell
_n^2}.  \label{ab}
\end{equation}
It is easy to verify
\[
\sum_kE_{k-1}f\otimes d_kg=D(\sum_{i<j\leq n}\alpha _i\beta _je_i\otimes
e_j)D.
\]
and
\begin{equation}
\left\| \sum_kE_{k-1}f\otimes d_kg\right\| _{L^1(S_n^1)}=\left\|
\sum_{i<j\leq n}\alpha _i\beta _je_i\otimes e_j\right\| _{S_n^1}=\left\|
T(\alpha \otimes \beta )\right\| _{S_n^1}.  \label{fgt}
\end{equation}
On the other hand, by duality between $L^1(S_n^1)$ and $L^\infty (M_n),$ we
have,
\begin{eqnarray}
\left\| \sum_kE_{k-1}f\otimes d_kg\right\| _{L^1(S_n^1)} &=&\sup \{\
tr\int \sum_kd_kb(E_{k-1}f\otimes d_kg), \Vert b\Vert _{L^\infty
(M_n)}\leq 1\}\nonumber \\
&\leq &\sup \{\ \Vert \pi _b(f)\Vert
_{L^2(\ell _n^2)}\Vert g\Vert _{L^2(\ell _n^2)}, \Vert b\Vert _{L^{^\infty }(M_n)}\leq 1\}  \nonumber \\
&\leq &c(n)\left\| f\right\| _{L^2(\ell _n^2)}\left\| g\right\|
_{L^2(\ell _n^2)}.  \label{fg}
\end{eqnarray}
Combining (\ref{fg}), (\ref{ab}) and (\ref{fgt}) we get (\ref{T}) and the
proof is complete. \qed

Recall that the square function of $b$ is defined as
\[
S(b)=(\sum_k|d_kb|^2)^{\frac 12}.
\]
The so called ``sweep'' function is just the square of the square function,
for this reason we denote it by $S^2(b)$,
\[
S^2(b)=\sum_k|d_kb|^2.
\]
In the classical case, we know that
\begin{eqnarray}
||S(b)||_{BMO_d} &\leq &c||b||_{BMO_d}  \label{sb} \\
||S^2(b)||_{BMO_d} &\leq &c||b||_{BMO_d}^2  \label{swb}
\end{eqnarray}
When considering square functions $S(b)$ for $M_n$ valued functions $b$,
a similar result remains true with an absolute constant.

\begin{proposition}
For any $n\in {\Bbb N}$, and any $M_n$ valued function $b$, we have
\[
||S(b)||_{BMO_c}\leq \sqrt{2}||b||_{BMO_c}
\]
\end{proposition}

\noindent {\bf Proof.} Since we are in the dyadic case, we have
\begin{eqnarray*}
||S(b)||^2_{BMO_c} &\leq&
2\sup_m||E_m[(S(b)-E_mS(b))^*(S(b)-E_mS(b))]||_{L^\infty(M_n)} \\
&=&2\sup_m||E_mS^2(b)-(E_mS(b))^2||_{L^\infty(M_n)}
\end{eqnarray*}
Note
\[
E_mS^2(b)-\sum_{k=1}^m|d_kb|^2\geq E_mS^2(b)-(E_mS(b))^2\geq 0 .
\]
We get
\begin{eqnarray*}
||S(b)||^2_{BMO_c} &\leq&
2\sup_m||E_mS^2(b)-\sum_{k=1}^m|d_kb|^2||_{L^\infty(M_n)} \\
&=&2\sup_m||E_m\sum_{k=m+1}|d_kb|^2||_{L^\infty(M_n)} \\
&\leq& 2||b||_{BMO_c}^2. \qed
\end{eqnarray*}

Matrix valued sweep functions have been studied in \cite{BS}, \cite{[GPTV]}
etc. Unlike in the case of square functions, it is proved in \cite{BS} that the best
constant $c_n$ such that
\begin{eqnarray}
||S^2(b)||_{BMO_c}\leq c_n||b||_{BMO_c}^2  \label{swm}
\end{eqnarray}
is $c\log (n+1)$. The following result shows that the best constant $c_n$ is
still $c\log (n+1)$ even if we replace $||\cdot ||_{BMO_c}$ by the bigger norm $%
||\cdot ||_{L^\infty (M_n)}$ in the right side of (\ref{swm}).

\begin{theorem}
For every $n\in {\Bbb N},$ there exists an $M_n$ valued function $b$ with $%
\left\| b\right\| _{L^\infty (M_n)}\leq 1$ but such that
\[
\left\| S^2(b)\right\| _{BMO_c}\geq c\log (n+1).
\]
\end{theorem}

\noindent{\bf Proof. }Consider a function $b$ that works for the statement of Theorem
1.1. Then $\Vert b\Vert _{L^\infty (M_n)}\leq 1$ and there exists a function $f\in
L^2(S_n^2),$ such that $\left\| f\right\| _{L^2(S_n^2)}\leq 1$ and
\begin{equation}
\left\| \sum_kd_kbE_{k-1}f\right\| _{L^2(S_n^2)}\geq c\log (n+1).
\label{log}
\end{equation}
We compute the square of the left side of (\ref{log}) and get
\begin{eqnarray*}
&&\left\| \sum_kd_kbE_{k-1}f\right\| _{L^2(S_n^2)}^2 \\
&=&tr\int \sum_k|d_kb|^2E_{k-1}fE_{k-1}f^{*} \\
&=&tr\int \sum_k|d_kb|^2(\sum_{i< k}|d_if^{*}|^2+\sum_{i<
k}E_{i-1}fd_if^{*}+\sum_{i< k}d_ifE_{i-1}f^{*}) \\
&=&tr\int \sum_i(\sum_{k> i}|d_kb|^2)|d_if^{*}|^2+tr\int
\sum_i(\sum_{k> i}|d_kb|^2)(E_{i-1}fd_if^{*}+d_ifE_{i-1}f^{*}) \\
&=&I+II
\end{eqnarray*}
For $I,$ note $|d_if^{*}|^2$ is ${\cal F}_i$ measurable, we have
\begin{eqnarray*}
I &=&tr\int \sum_iE_i(\sum_{k> i}|d_kb|^2)|d_if^{*}|^2 \\
&\leq &\sup_i||E_i(\sum_{k> i}|d_kb|^2)||_{L^\infty (M_n)} (tr\int \sum_i
|d_if^{*}|^2) \\
&\leq &||b||_{BMO_c}^2||f||_{L^2(S_n^2)}^2\leq 4
\end{eqnarray*}
For $II,$ note $E_{i-1}fd_if^{*}+d_ifE_{i-1}f^{*}$ is a martingale
difference and $\sum_{k\leq i}|d_k|^2$ is ${\cal F}_{i-1}$ measurable since we are in the dyadic case,
we get
\begin{eqnarray*}
II &=&tr\int \sum_iS^2(b)(E_{i-1}fd_if^{*}+d_ifE_{i-1}f^{*}) \\
&=&tr\int \sum_id_i(S^2(b))(E_{i-1}fd_if^{*}+d_ifE_{i-1}f^{*}) \\
&\leq &2||\sum_id_i(S^2(b))E_{i-1}f||_{L^2(S_n^2)}||f||_{L^2(S_n^2)} \\
&\leq &2||\pi _{S^2(b)}||_{L^2(S_n^2)\rightarrow L^2(S_n^2)} \\
&\leq &2c\log (n+1)||S^2(b)||_{BMO_c}.
\end{eqnarray*}
We used (\ref{katz}) in the last step. Combining this with (\ref{log}), we
get
\[
c\log (n+1)\leq \left\| \sum_kd_kbE_{k-1}f\right\| _{L^2(S_n^2)}^2\leq
4+2c\log (n+1)||S^2(b)||_{BMO_c}
\]
Thus
\[
||S^2(b)||_{{\rm BMO}_c}\geq c\log (n+1).
\]
This completes the proof. \qed

\section{Proof of Theorem 1.2.}

We keep the notations introduced in the end of Section 1. Recall BMO spaces
of noncommutative martingales are defined for $x=(x_k)\in L^2({\cal M})$ as
below (see \cite{PX}, \cite{JX}):
\begin{eqnarray*}
{\rm BMO}_c({\cal M}) &=&\{x:||x||_{{\rm BMO}_c({\cal M})}=\sup_n\left\|
E_n|\sum_{k=n}^\infty d_kx|^2\right\| _{{\cal M}}^{\frac 12}<\infty \}; \\
{\rm BMO}_r({\cal M}) &=&\{x:||x||_{{\rm BMO}_r({\cal M})}=||x^{*}||_{{\rm %
BMO}_c({\cal M})}<\infty \}; \\
{\rm BMO}_{cr}({\cal M}) &=&\{x:||x||_{{\rm BMO}_{cr}({\cal M})}=\max
\{||x||_{{\rm BMO}_c({\cal M})},||x||_{{\rm BMO}_r({\cal M})}\}<\infty \}.
\end{eqnarray*}

When ${\cal M}=L^\infty (M_n),$ ${\rm BMO}_c({\cal M})$ is just ${\rm BMO}_c$
considered in Section 1 and 2. In this section, for noncommutative
martingale $b,$ we consider $\pi _b$ and $\widetilde{\pi }_b$ as operators
on bounded noncommutative $L^p$-martingale spaces introduced in Section 1.
We will need the following interpolation result and the John-Nirenberg theorem
for noncommutative martingales proved by Junge and Musat recently (see \cite
{JM}, \cite{MM}).

\begin{theorem}
(Musat) For $1\leq p\leq q<\infty ,$
\[
(BMO_{cr}({\cal M}),L_p({\cal M}))_\theta =L_q({\cal M}),\mbox{ with }\theta
=\frac{p}q.
\]
\end{theorem}

\begin{theorem}
(Junge, Musat) For any $1\leq q<\infty$ and any $g=(g_k)_k\in
BMO_{cr}({\cal M})$, there exist $c_q, c'_q>0$ such that
\begin{equation}
c_q^{\prime }||g||_{BMO_{cr}}\leq \sup_{m\in {\Bbb N}}\sup_{a\in
{\cal M}_m, \tau(|a|^q)\leq 1}\{||\sum_{k\geq m}d_kga||_{L^q({\cal
M})},||\sum_{k\geq m}ad_kg||_{L^q({\cal M})}\}\leq
c_q||g||_{BMO_{cr}}.  \label{jn}
\end{equation}
\end{theorem}
In fact, the formula above is proved for $q\geq 2$ in \cite{JM}. It is not
hard to show that it is also true for $1\leq q<2.$ In the following, we
give a simpler proof of it in the tracial case.

\noindent {\bf Proof.} Note for any $g\in BMO_{cr}({\cal M})$,
\[
||g||_{BMO_{cr}({\cal M})}=\sup_{m\in {\Bbb N}} \sup_{a\in {\cal
M}_m, \tau(|a|^2)\leq 1}\{||\sum_{k\geq m}d_kga||_{L^2({\cal
M})},||\sum_{k\geq m}ad_kg||_{L^2({\cal M})}\}.
\]
We get $c_2=c_2^{\prime }=1.$ Note for $p,r,s$ with $1/p=1/r+1/s$ and
$a\in L^p({\cal M}), ||a||_{L^p({\cal M})}\leq 1$, there exist
$b, c$ such that
$a=bc$ and $||b||_{L^r({\cal M})}\leq 1, ||c||_{L^s({\cal M})}\leq 1$.
By H$\stackrel{..}{o}$lder's inequality we then get $c_q=1$ for $1\leq q<2$ and $c_q^{\prime }=1$ for $2<q<\infty .$ Thus
for $2<q<\infty ,$ we only need to prove the second inequality of (\ref{jn}%
). And, for $1\leq q<2$, we only need to prove the first inequality of (\ref
{jn}). Fix $g\in BMO_{cr}({\cal M}),$ $m\in {\Bbb N},$ consider the left
multiplier $L_m$ and the right multiplier $R_m$ defined as
\[
L_m(a)=\sum_{k\geq m}d_kga\mbox{ and }R_m(a)=\sum_{k\geq m}ad_kg,\quad
\forall a\in {\cal M}_m.
\]
It is easy to check that
\begin{eqnarray*}
\sup_m||L_m||_{L^2({\cal M}_m)\rightarrow L^2({\cal M})} &=&||g||_{BMO_c}, \\
\sup_m||L_m||_{L^\infty ({\cal M}_m)\rightarrow BMO_{cr}} &\leq
&||g||_{BMO_{cr}}; \\
\sup_m||R_m||_{L^2({\cal M}_m)\rightarrow L^2({\cal M})} &=&||g||_{BMO_r},
\\
\sup_m||R_m||_{L^\infty ({\cal M}_m)\rightarrow BMO_{cr}} &\leq
&||g||_{BMO_{cr}}.
\end{eqnarray*}
Thus $L_m,R_m$ extend to bounded operators from $L^2({\cal M}_m)$ to $L^2(%
{\cal M}),$ as well as from $L^\infty ({\cal M}_m)$ to
$BMO_{cr}({\cal M})$. By Musat's interpolation result Theorem 3.5,
we get $L_m$ and $R_m$ are bounded from $L^q({\cal M}_m)$ to
$L^q({\cal M})$ and their operator norms are smaller than
$c_q||g||_{BMO_{cr}},$ for all $2\leq q<\infty .$ By taking supremum
over $m,$ we prove the second inequality of (\ref{jn}) for $q\geq
2$.

For $1\leq q<2,$ by interpolation again, for $\theta =\frac{q}2$ and some $%
c_q^{\prime \prime }>0,$
\begin{eqnarray*}
||L_m||_{L^2({\cal M}_m)\rightarrow L^2({\cal M})} &\leq &c_q^{\prime \prime
}||L_m||_{L^q({\cal M}_m)\rightarrow L^q({\cal M})}^\theta ||L_m||_{L^\infty
({\cal M}_m)\rightarrow BMO_{cr}}^{1-\theta } \\
&\leq &c_q^{\prime \prime }||L_m||_{L^q({\cal M}_m)\rightarrow L^q({\cal M}%
)}^\theta ||g||_{BMO_{cr}}^{1-\theta }, \\
||R_m||_{L^2({\cal M}_m)\rightarrow L^2({\cal M})} &\leq &c_q^{\prime \prime
}||R_m||_{L^q({\cal M}_m)\rightarrow L^q({\cal M})}^\theta ||R_m||_{L^\infty
({\cal M}_m)\rightarrow BMO_{cr}}^{1-\theta } \\
&\leq &c_q^{\prime \prime }||R_m||_{L^q({\cal M}_m)\rightarrow L^q({\cal M}%
)}^\theta ||g||_{BMO_{cr}}^{1-\theta }.
\end{eqnarray*}
Thus
\begin{eqnarray*}
||g||_{BMO_{cr}} &=&\max \{\sup_m||L_m||_{L^2({\cal M}_m)\rightarrow L^2(%
{\cal M})},\sup_m||R_m||_{L^2({\cal M}_m)\rightarrow L^2({\cal M})}\} \\
&\leq &c_q^{\prime \prime }||g||_{BMO_{cr}}^{1-\theta }\sup_m\{||L_m||_{L^q(%
{\cal M}_m)\rightarrow L^q({\cal M})}^\theta ,||R_m||_{L^q({\cal M}%
_m)\rightarrow L^q({\cal M})}^\theta \}.
\end{eqnarray*}
This gives the first inequality of (\ref{jn}) with $c_q^{\prime
}=(c_q^{\prime \prime })^{-\frac 1\theta }$ for $1\leq q<2.$ \qed

Recall that we say a filtration ${\cal M}_k$ is ``regular" if, for
some $c>0$, $||a||_\infty\leq c||E_{m-}a||_\infty, \ \forall m\in
{\Bbb N}, a\geq 0, a\in {\cal M}_m$.
\begin{lemma}
For any regular filtration ${\cal M}_k$, we have
\begin{equation}
||b||_{BMO_{cr}({\cal M})}\leq c_p\max \{||\pi_b||_{L^p({\cal M}%
)\rightarrow L^p({\cal M})},||\widetilde{\pi }_b||_{L^p({\cal M})\rightarrow
L^p({\cal M})}\},\ \ \forall 1\leq p<\infty .  \label{last}
\end{equation}
\end{lemma}
\noindent
{\bf Proof.}
Note, for any $b\in BMO_{cr}({\cal M})$ with respect to the regular filtration ${\cal M}_k$,
\begin{eqnarray*}
||b||_{BMO_{cr}({\cal M})}\leq  c\sup_{m\in {\Bbb N}}\sup_{\tau
a^2\leq 1,a\in {\cal M}_m}\{||\sum_{k > m}d_kba||_{L^2({\cal
M})},||\sum_{k> m}ad_kb||_{L^2({\cal M})}\}.
\end{eqnarray*}
Similar to the proof of Theorem 3.6, we can get,
\begin{equation}
c_q^{\prime }||b||_{BMO_{cr}}\leq \sup_{m\in {\Bbb N}} \sup_{a\in
{\cal M}_m, \tau |a|^q\leq 1} \{||\sum_{k>m}d_kba||_{L^q({\cal
M})},||\sum_{k> m}ad_kb||_{L^q({\cal M})}\}\leq c_q||b||_{BMO_{cr}}.
\label{jnd}
\end{equation}
On the other hand, by considering $\pi_b(a), \widetilde{\pi}_b(a)$
for $a\in {\cal M}_m$, $||a||_{L^p({\cal M})}\leq 1$, we have
\begin{eqnarray*}
\sup_{a\in {\cal M}_m, \tau |a|^q\leq 1}\{||\sum_{k>
m}d_kba||_{L^p({\cal M})},||\sum_{k>
m}ad_kb||_{L^p({\cal M})}\}\\
\leq 2\max \{||\pi_b||_{L^p({\cal M}%
)\rightarrow L^p({\cal M})},||\widetilde{\pi }_b||_{L^p({\cal M})\rightarrow
L^p({\cal M})}\}.
\end{eqnarray*}
Taking supremum over $m$ in the inequality above, we get (\ref{last}) by (\ref{jnd}) . \qed

\begin{lemma}
For $1<p<\infty ,$ we have
\begin{eqnarray}
\left\| \pi _b\right\| _{L^\infty ({\cal M})\rightarrow BMO_{cr}({\cal M}
)}&\leq& c_p(\left\| \pi _b\right\| _{L^p({\cal M})\rightarrow L^p({\cal M}
)}+||b||_{BMO_r({\cal M})}). \label{1}\\
\left\| \widetilde \pi _b\right\| _{L^\infty ({\cal M})\rightarrow BMO_{cr}({\cal M}
)}&\leq& c_p(\left\| \widetilde \pi _b\right\| _{L^p({\cal M})\rightarrow L^p({\cal M})}
+||b||_{BMO_c({\cal M})}). \label{2}
\end{eqnarray}
\end{lemma}

\noindent{\bf Proof. } We prove (\ref{1}) only. Fix a $f\in L^\infty ({\cal M})$ with $\left\|
f\right\| _{L^\infty ({\cal M})}\leq 1.$ We have
\begin{eqnarray}
&&\left\| E_m\sum_{k\geq m}|d_kbE_{k-1}f|^2\right\| _{L^\infty ({\cal M})}
\nonumber \\
&=&\sup \{\tau E_m\sum_{k\geq m}|d_kbE_{k-1}f|^2a,\ a\in {\cal M}_m,a\geq
0,\tau a\leq 1\}  \nonumber \\
&=&\sup\{\tau \sum_{k\geq
m}(d_kbE_{k-1}fa^{\frac 1p})^{*}(d_kbE_{k-1}fa^{\frac 1q}),\ a\in {\cal M}_m,a\geq
0,\tau a\leq 1\}   \nonumber \\
&\leq &\sup_a\left\| d_mbE_{m-1}fa^{\frac
1p}+\sum_{k>m}d_kbE_{k-1}(fa^{\frac 1p})\right\| _{L^p({\cal M})}\left\|
\sum_{k\geq m}d_kbE_{k-1}fa^{\frac 1q}\right\| _{L^q({\cal M}{\Bbb {)}}}
\nonumber
\end{eqnarray}
Note $||d_mbE_{m-1}fa^{\frac1p}||_{L^p({\cal M})}\leq ||d_mb||_{\cal M}
\leq {||b||_{BMO_r}}$.
By (\ref{jn}) we get
\begin{eqnarray}
\left\| E_m\sum_{k\geq m}|d_kbE_{k-1}f|^2\right\| _{L^\infty ({\cal M})}
\leq c_q(||b||_{BMO_r}+\left\| \pi _b\right\| _{L^p({\cal M})\rightarrow L^p({\cal M}%
)})\left\| \pi _b(f)\right\| _{BMO_{cr}({\cal M)}}.  \label{bmoc}
\end{eqnarray}
Taking supremum over
$m$ in (\ref{bmoc}), we get
\[
\left\| \pi _b(f)\right\| _{BMO_c({\cal M})}^2
\leq c_q(||b||_{BMO_r}+\left\| \pi _b\right\| _{L^p({\cal M})\rightarrow L^p({\cal M}%
)})\left\| \pi _b(f)\right\| _{BMO_{cr}({\cal M)}}.
\]
On the other hand, since $(E_{m-1}f)(E_{m-1}f)^*\leq 1$, we have
\[
\left\| \pi _b(f)\right\| _{BMO_r({\cal M})}\leq \left\| b\right\| _{BMO_r(%
{\cal M})}.
\]
Thus,
\[
\left\| \pi _b(f)\right\| _{BMO_{cr}({\cal M})}^2\leq (c_q+1)(\left\| \pi
_b\right\| _{L^p({\cal M})\rightarrow L^p({\cal M})}+||b||_{BMO_r({\cal M}%
)})\left\| \pi _b(f)\right\| _{BMO_{cr}({\cal M})},%
\]
Therefore
\[
\left\| \pi _b\right\| _{L^\infty ({\cal M})\rightarrow BMO_{cr}({\cal M}%
)}\leq (c_q+1)(\left\| \pi _b\right\| _{L^p({\cal M})\rightarrow L^p({\cal M}%
)}+||b||_{BMO_r({\cal M})}).\qed
\]

\medskip {\it Proof of Theorem 1.2.} By Lemma 3.7 and Lemma 3.8 we get immediately that
\begin{eqnarray*}
\max&\{\left\| \pi _b\right\| _{L^\infty({\cal M})\rightarrow BMO_{cr}},
\left\| \widetilde \pi _b\right\| _{L^\infty({\cal M})\rightarrow BMO_{cr}}\}\\
\leq
c_p\max&\{\left\| \pi _b\right\| _{L^p({\cal M})\rightarrow L^p({\cal M}%
)}, \ \ \left\| \widetilde \pi _b\right\| _{L^p({\cal M})\rightarrow L^p({\cal M}%
)}\}
\end{eqnarray*}
By the interpolation results of
noncommutative martingales( Theorem 3.5), we get
\begin{eqnarray*}
\max&\{\left\| \pi _b\right\| _{L^q({\cal M})\rightarrow L^q({\cal M})},
\left\| \widetilde\pi _b\right\| _{L^q({\cal M})\rightarrow L^q({\cal M})}\}\\
\leq
c_p \max&\{\left\| \pi _b\right\| _{L^p({\cal M})\rightarrow L^p({\cal M}%
)}, \left\| \widetilde\pi _b\right\| _{L^p({\cal M})\rightarrow L^p({\cal M}%
)}\},
\end{eqnarray*}
for all $1<p<q<\infty.$

\medskip

 {\bf Question :} Assume ${\pi }_b,
\widetilde{\pi }_b$ are of type $(p,p)$, are they of weak type
$(1,1)$? More precisely, assume $||\pi_b||_{L^p({\cal M}
)\rightarrow L^p({\cal M})}+||\widetilde{\pi }_b||_{L^p({\cal
M})\rightarrow L^p({\cal M})}< \infty$, does there exist a constant
$C>0$ such that, for any $f\in L^1({\cal M})$, $\lambda>0$, there is
a projection $e\in {\cal M}$ such that
$$
\tau(e^{\bot})\leq C\frac {||f||_{L^1({\cal M})}}{\lambda}\ \ \
{\mathrm and}\ \ \ ||e{\pi }_b(f)e||_{L^\infty({\cal
M})}+||e\widetilde{\pi }_b(f)e||_{L^\infty({\cal M})}\leq \lambda ?
$$
\medskip

 We have the following corollary by applying results of this
section to matrix valued dyadic paraproducts discussed in Section 1
and Section 2. Note $M_n$ valued dyadic martingales on the unit
circle are noncommutative martingales associated with the von Neuman
algebra ${\mathcal M}=L^\infty({\Bbb T})\otimes M_n$ and the
filtration ${\mathcal M}_k=L^\infty({\Bbb T}, {\mathcal F}_k)\otimes
M_n$.

\begin{corollary}
Let $1<p<\infty ,$ denote by $c_p(n)$ the best constant such that
\[
\left\| \pi _b\right\| _{L^p(S_n^p)\rightarrow L^p(S_n^p)}\leq c_p(n)\left\|
b\right\| _{L^\infty (M_n)},\ \forall b.
\]
Then
\[
c_p(n)\backsim \log (n+1).
\]
\end{corollary}

\noindent{\bf Proof. }Note in the proof of Theorem 1.1, if we see $f$
as a column matrix valued function and $g$ as a row matrix valued function,
we will have
\[
||f||_{L^p(S_n^p)}=||\alpha ||_{\ell _n^2},\ \ ||g||_{L^q(S_n^q)}=||\beta
||_{\ell _n^2}.
\]
By the same method, we can prove $c_p(n)\geq c\log (n+1)$ for all $1<p<\infty$.
For the inverse relation, by (\ref{katz}) we have $c_2(n)\leq c\log(n+1)$. Then, by (\ref{1}),
we get
\begin{eqnarray}
\left\| \pi _b\right\| _{L^\infty (M_n)\rightarrow BMO_{cr}}
&\leq&
c_2(c_2(n)\left\| b\right\| _{L^\infty (M_n)}+||b||_{BMO_{cr}}) \nonumber\\
&\leq& c\log(n+1)||b||_{L^\infty(M_n)}, \ \ \ \forall b\in L^\infty(M_n) \label{b}
\end{eqnarray}
Denote by $\pi^*_b$ the adjoint operator of the dyadic paraproduct $\pi_b$, then
$$
\pi_b^*(f)=\sum_k(d_kb)^*E_{k-1}f.
$$
Note we have the decomposition
\[
\pi _b^{*}(f)=b^{*}f-\pi _{b^{*}}(f)-(\pi _{f^{*}}(b))^{*}.
\]
By (\ref{b}), we get
\begin{eqnarray}
\left\| \pi _b^{*}\right\| _{L^\infty (M_n)\rightarrow BMO_{cr}}
&\leq& ||b^*||_{L^\infty (M_n)}+c\log(n+1)||b^*||_{L^\infty (M_n)}+c\log(n+1)||b||_{L^\infty (M_n)}
\nonumber\\
&\leq&
c\log(n+1)\left\| b\right\| _{L^\infty (M_n)}. \label{b*}
\end{eqnarray}
By (\ref{b}), (\ref{b*}) and the interpolation result Theorem 3.5, we get
\[
\left\| \pi _b\right\| _{L^p(S_n^p)\rightarrow L^p(S_n^p)}\leq
c_p\log(n+1)\left\| b\right\| _{L^\infty (M_n)},\ \ \ \forall 1<p<\infty.
\]
Therefore, we can conclude $c_p(n)\backsim
\log (n+1).$ \qed

\vskip 0.2 in \noindent {\bf Acknowledgments.} The author is
indebted to his adviser G. Pisier for many helpful conversations.
The example in the proof of Theorem 1.1 comes out from the
discussions with him. The author also would like to thank the
referee for his very useful comments.

\medskip
$
\begin{array}{l}
\mbox{Dept. of Math.} \\
\mbox{Texas A\&M University} \\
\mbox{College Station, TX, 77843} \\
\mbox{U. S. A.} \\
\mbox{tmei@math.tamu.edu}
\end{array}
$

\end{document}